\documentclass[a4paper,reqno,12pt]{amsart}   
\usepackage{newtxtext} 
\usepackage{amsmath}
\usepackage[bigdelims]{newtxmath}
\usepackage[T1]{fontenc}
\usepackage{textcomp}

\usepackage{tikz}
\usepackage[mathscr]{eucal}

\newcounter{mylist}

\numberwithin{equation}{section}	

\theoremstyle{plain}

\newtheorem{theorem}{Theorem}[section]      
\newtheorem*{theorem*}{Theorem}

\newtheorem{lemma}[theorem]{Lemma}

\newtheorem{proposition}[theorem]{Proposition}

\newtheorem{corollary}[theorem]{Corollary}

\newtheorem*{conjecture*}{Conjecture}

\theoremstyle{definition}
\newtheorem{definition}[theorem]{Definition}

\newtheorem{example}[theorem]{Example}

\theoremstyle{remark}
\newtheorem{remark}[theorem]{Remark}

\def\geqsl{\geqslant}
\def\leqsl{\leqslant}

\renewcommand{\setminus}{\smallsetminus}  

\newcommand{\trace}[1]{\ensuremath{\operatorname{tr}{\left( #1 \right)}}}

\newcommand{\tp}[1]{\ensuremath{\hspace{1.8truept}\vphantom{1}^{t}\hspace{-1pt}#1}}

\renewcommand{\Xi}{\ensuremath{\varXi}}
\renewcommand{\Theta}{\ensuremath{\varTheta}}

\newcommand{\N}{\ensuremath{\mathbb N}}
\newcommand{\Z}{\ensuremath{\mathbb Z}}
\newcommand{\Q}{\ensuremath{\mathbb Q}}
\newcommand{\R}{\ensuremath{\mathbb R}}
\newcommand{\C}{\ensuremath{\mathbb C}}

\newcommand{\mathfraksl}{\ensuremath{\mathfrak {sl}}}

\newcommand{\g}{\ensuremath{\mathfrak{g}}}   

\newcommand{\GL}[1]{\ensuremath{\mathrm{GL}_{#1}}}  
\newcommand{\mathO}{\ensuremath{\mathrm{O}}}          
\newcommand{\U}{\ensuremath{\mathrm{U}}}            
\newcommand{\Sp}{\ensuremath{\mathrm{Sp}}}          
\newcommand{\Mp}{\ensuremath{\mathrm{Mp}}}          

\def\Z{\mathbb Z}

\newcommand{\rmd}{\ensuremath{\operatorname{d}\!}}
\newcommand{\pd}{\ensuremath{\partial}}

\def\<{\langle}
\def\>{\rangle}

\def\c2vec#1#2{ %
   \left[ \begin{smallmatrix} %
           #1 \\ #2  \end{smallmatrix} %
   \right]}

\newcommand{\ibar}{\ensuremath{\bar{\imath}}}
\newcommand{\jbar}{\ensuremath{\bar{\jmath}}}
\DeclareMathOperator{\impart}{Im}
\DeclareMathOperator{\repart}{Re}

\DeclareMathOperator{\spanned}{span}
\DeclareMathOperator{\GKdim}{Dim}
\DeclareMathOperator{\B-deg}{Deg}

\newcommand{\rmO}{\mathO}
\newcommand{\ai}{\ensuremath{\mathrm{i}}\,}
\newcommand{\rfac}[2]{\ensuremath{(#1)_{#2}}}
\newcommand{\ffac}[2]{\ensuremath{(#1)^{-}_{#2}}}

\newcommand{\abs}[1]{\ensuremath{\left|\, #1 \, \right|}}

\newcommand{\PV}{\ensuremath{\mathscr P(V)}}

\newcommand{\PD}{\ensuremath{\mathscr{PD}}}

\newcommand{\Eu}{\ensuremath{E}}

\newcommand{\Harm}{\ensuremath{\mathscr H}}
\newcommand{\ourE}{\ensuremath{\mathscr E}}
\newcommand{\Module}[2]{\ensuremath{\mathit M}^{#1}(#2)}
\newcommand{\suppKtypeModule}{\ensuremath{\Sigma}}

\newcommand{\Cspan}[1]{\mbox{\C -}\! \spanned \left\{ {#1} \right\} }

\begin{document}

\title[$(\g,K)$-module of $\rmO(p,q)$]
{$(\g,K)$-module of $\rmO(p,q)$ associated with the finite-dimensional representation of $\mathfraksl_2$}

\author{Takashi Hashimoto}
\thanks{Partly supported by JSPS Grant-in-Aid for Scientific Research (C) No.~23540203 and No.~26400014.}

\address{
  University Education Center, 
  Tottori University, 
  4-101, Koyama-Minami, Tottori, 680-8550, Japan
}
\email{thashi@tottori-u.ac.jp}
\date{\today}

\begin{abstract}
The main aim of this paper is 
to construct irreducible $(\g,K)$-modules of $\rmO(p,q)$ corresponding to the finite-dimensional representation of $\mathfraksl_2$ of dimension $m+1$ under the Howe duality,
to find the $K$-type formula, the Gelfand-Kirillov dimension and the Bernstein degree of them for non-negative integers $m$. 
The $K$-type formula for $m=0$ shows that it is nothing but the $(\g,K)$-module of the minimal representation of $\rmO(p,q)$.
One finds that the Gelfand-Kirillov dimension is equal to $p+q-3$ not only for $m=0$ but for any $m$ satisfying $m+3 \leqsl (p+q)/2$ when $p, q \geqsl 2$ and $p+q$ is even,
and that the Bernstein degree for $m$ is equal to $(m+1)$ times that for $m=0$.
\end{abstract}

\maketitle

\section{Introduction}

Let $G$ be a Lie group with $\g_0$ its Lie algebra and $\g$ the complexification of $\g_0$.
An action of $G$ on a symplectic manifold $(M,\omega)$ is called symplectic if $g^* \omega = \omega$ for all $g \in G$,
and a symplectic action is called Hamiltonian if there exists a smooth $G$-equivariant map $\mu: M \to \g_0^*$ satisfying the condition \eqref{e:definition of moment map} below, 
which is called a moment map, where $\g_0^*$ is the dual vector space of $\g_0$.
We are concerned with the cases where the symplectic manifold is a real symplectic vector space $(W,\omega)$.
It was shown in \cite{Harm_oscil_kjm} that when $G=\Sp(n,\R), \U(p,q)$ and $\rmO^*(2n)$, 
the canonical quantization of the moment map on $W=\R^{2n}, (\C^{p+q})_{\R}$ and $(\C^{2n})_{\R}$, with a choice of a Lagrangian subspace in each case, yields a representation of $\g$ 
that is the differentiation of the oscillator (or Segal-Shale-Weil) representation of $\Mp(n,\R), \U(p,q)$ and $\rmO^*(2n)$ respectively, 
where $\Mp(n,\R)$ is the metaplectic group, i.e.,~the double cover of $\Sp(n,\R)$.

In this paper, we consider the case where $W=(\C^{p+q})_\R$, the real vector space underlying $\C^{p+q}$: 
\[
  W = \left\{ z=x + \ai y \left|\, x, y \in \R^{p+q} \right. \right\},
\]
which we regard as a symplectic vector space equipped with a symplectic form $\omega$ given by 
\[
  \omega(z,w)=\impart( z^* I_{p,q} w)  \qquad (z, w \in W),
\]
and $G=\rmO(p,q)$, the indefinite orthogonal group defined by
\[
  \rmO(p,q) = \left\{ g \in \GL{p+q}(\R) \left| \tp{g} I_{p,q} g = I_{p,q} \right. \right\}
\]
with $I_{p,q} = \left[ \begin{smallmatrix} 1_p & \\  & -1_q \end{smallmatrix} \right]$.
The action of $G=\rmO(p,q)$ on $W$ defined by matrix multiplication is symplectic and Hamiltonian.
The $\rmO(p,q)$-case we consider here is closely related to the $\U(p,q)$-case mentioned above.
In fact, the symplectic vector space $(W,\omega)$ for $\rmO(p,q)$ is identical to the one for $\U(p,q)$, 
and the action of $\rmO(p,q)$ on $W$ is the restriction of the action of $\U(p,q)$ induced from the canonical embedding of $\rmO(p,q)$ into $\U(p,q)$.
Furthermore, the moment map for the $\rmO(p,q)$-case is the real part of the one for the $\U(p,q)$-case.

The canonical quantization of the moment map $\mu: W=(\C^{p+q})_\R \to \g_0^*$ for $G=\rmO(p,q)$, with a choice of a Lagrangian subspace $V$ of $W$, 
yields a representation $\pi$ of $\g$ as in the cases mentioned above, which is shown to be a partial Fourier transformation of the representation $\pi^\sharp$ of $\g$
obtained by differentiating the left regular representation of $G$ on $C^{\infty}(V)$.
Note that if we restrict the operator $\pi^\sharp(X)$, $X \in \g$, to a subspace consisting of homogeneous functions on $V$ with respect to the multiplicative group $\R_{>0}$, 
then the restricted representation is the degenerate principal series of $G$ obtained by inducing up a one-dimensional representation of a parabolic subgroup of $G$ (see \cite{HT93}).

In the influential paper \cite{Howe_remarks}, Howe showed that one can treat the classical invariant theory from a unified viewpoint --- the dual pair.
We focus our attention on the dual pair $(\rmO(p,q), \mathrm{SL}_2(\R))$, both components of which are non-compact, 
and apply the representation theory of $\mathfrak{sl}_2$ to cut out irreducible $(\g,K)$-modules, which we denote by $\Module{+}{m}$ and $\Module{-}{m}$, $m=0,1,2,\dots$, in this paper,
where $\Module{+}{m}$ (resp.~$\Module{-}{m}$) consists of all highest (resp.~lowest) weight vectors with respect to the $\mathfrak{sl}_2$-action 
(see Definition \ref{d:Module^{pm}_{m}} below for details).
We will see that such weight vectors are given in terms of harmonic polynomials and the Bessel functions of the first kind.
Both $\Module{\pm}{m}$ correspond to the $(m+1)$-dimensional irreducible representation of $\mathfrak{sl}_2$ under the Howe duality, and in fact are isomorphic to each other.
They were originally considered in \cite{RS80} without the condition of finite-dimensionality.
Note that $\Module{+}{0} = \Module{-}{0}$ by definition.

In the cases of the oscillator representations mentioned above, i.e.,~when $G=\Sp(n,\R)$, $\U(p,q)$ and $\rmO^*(2n)$,
we note that the counterpart $G'$ of $G$ for the dual pair $(G,G')$ is compact, hence, all its irreducible representations are finite-dimensional.
Furthermore, the oscillator representations give examples of the minimal representations (we refer to \cite{KM2011} and the references therein for the definition of the minimal representation).
When $G=\rmO(p,q)$, its minimal representation is discussed e.g.~in \cite{Kostant90, BZ91, ZH97, Kobayashi_Orsted_03-1-3, KM2011}.

The main result of this paper is the $K$-type formula of $\Module{\pm}{m}$ for non-negative integers $m$ satisfying 
\begin{equation}
\label{e:condition_on_m}
  m+3 \leqsl \frac{p+q}{2},  
\end{equation} 
from which one can show that $\Module{\pm}{m}$ are irreducible $(\g,K)$-modules for $p, q \geqsl 2$ with $p+q$ even (Theorem \ref{t:main_result}).
The fact that the elements of $\Module{\pm}{m}$ are described in terms of the Bessel function plays a r\^{o}le in the proof of our main result.
The $K$-type formula of $\Module{+}{0}=\Module{-}{0}$, which corresponds to the one-dimensional trivial representation of $\mathfrak{sl}_2$, shows that 
it is nothing but the underlying $(\g,K)$-module of the minimal representation of $\rmO(p,q)$.
We will see that the Gelfand-Kirillov dimension of $\Module{\pm}{m}$ is equal to $p+q-3$ not only for $m=0$ but for {\em any} non-negative integer $m$ satisfying \eqref{e:condition_on_m}.
Meanwhile, the Bernstein degree of $\Module{\pm}{m}$ is $(m+1)$ times that of $\Module{+}{0}=\Module{-}{0}$ (Corollary \ref{c:GKdim_and_B-deg}).

The rest of this paper is organized as follows.
In \S\ref{s:moment_map_quantization}, 
we compute the moment map $\mu$ on $W$ for $G=\rmO(p,q)$, and construct the representation $\pi$ of $\g$ via canonical quantization of $\mu$.
Then we show that $\pi$ is a partial Fourier transform of the differential representation of the left regular representation of $G$ on $C^\infty(V)$.
In \S\ref{s:dual_pair_in_this_paper}, 
we give an $\mathfrak{sl}_2$-action that commutes with $\pi$, and find both highest weight vectors and lowest weight vectors with respect to the $\mathfrak{sl}_2$-action.
We remark that such weight vectors are given in terms of the Bessel functions of the first kind.
In \S\ref{s:(g,K)-module}, 
we introduce $(\g,K)$-modules $\Module{\pm}{m}$ and prove that $\Module{+}{m}$ and $\Module{-}{m}$ are isomorphic to each other for any non-negative integer $m$.
Then we find the $K$-type formula of $\Module{\pm}{m}$ for $m$ satisfying \eqref{e:condition_on_m} and show that they are irreducible.
As a corollary, we obtain the Gelfand-Kirillov dimension and the Bernstein degree of $\Module{\pm}{m}$.

\subsection*{Notation}
 
Let $\N$ denote the set of non-negative integers $\{0,1,2,\dots \}$, and $[p]$ the set $\{1,2,\dots,p\}$.
For the sake of simplicity, we write $\ibar := p+i$ for $i \in [q]$.
Finally, for $\alpha \in \C$ and $n \in \N$, we denote the rising and the falling factorials by
\[
  \rfac{\alpha}{n} := \prod_{i=1}^n (\alpha + i-1 )
     \quad  \text{and}  \quad
  \ffac{\alpha}{n} := \prod_{i=1}^n (\alpha - i + 1),
\]
respectively.

\section{Moment Map and its Quantization}
\label{s:moment_map_quantization}


Let $G$ be the indefinite orthogonal group $\rmO(p,q)$, which we realize by
\[
  \rmO(p,q) = \left\{ g \in \GL{p+q}(\R) \left| \tp{g} I_{p,q} g = I_{p,q} \right. \right\}
\]
with $I_{p,q}=\left[\begin{smallmatrix} 1_p &  \\  & -1_q \end{smallmatrix}\right]$.
Let $K$ be a maximal compact subgroup of $G$ given by 
\[
  K = \left\{ \left. \begin{bmatrix} a & 0 \\ 0 & d \end{bmatrix} \in G \, \right| a \in \rmO(p), d \in \rmO(q) \right\}
    \simeq \rmO(p) \times \rmO(q).
\]
We denote the Lie algebra of $K$ and its complexification by ${\mathfrak k}_0$ and $\mathfrak k$ respectively.

Let $\{X^{\pm}_{i,j} \}$ be a basis for $\g_0=\mathfrak{o}(p,q)$ given by
\begin{equation}
\label{e:basis}
\begin{aligned}
  X^+_{i,j} &= E_{i,j} - E_{j,i}                          & \quad  &(i,j \in [p]) \\
  X^+_{\ibar,\jbar} &= E_{\ibar,\jbar} - E_{\jbar,\ibar}  & \quad  &(i,j \in [q]) \\
  X^-_{i,j} &= E_{i,\jbar} + E_{\jbar,i}                  & \quad  &(i \in [p], j \in [q]),
\end{aligned}
\end{equation}
which also forms a basis for $\g=\mathfrak{o}_{p+q}$, the complexification of $\g_0=\mathfrak{o}(p,q)$.
We often identify $\g^*$ with $\g$ via the invariant bilinear form $B$ given by
\[
  B(X,Y) = \frac12 \trace{X Y}   \qquad (X,Y \in \g),
\]
where $\g^*$ denotes the dual space of $\g$.
Finally, let $\g=\mathfrak k \oplus \mathfrak p$ be the complexified Cartan decomposition of $\g$ with
\[
  \mathfrak k = \sum_{i,j \in [p]} \C X^+_{i,j} \oplus \sum_{i,j \in [q]} \C X^+_{\ibar,\jbar},
     \qquad
  \mathfrak p=\sum_{i \in [p], j \in [q]} \C X^-_{i,j}.
\]

Let $W$ be the real vector space $(\C^{p+q})_\R$ underlying the complex vector space $\C^{p+q}$:
\[
  W = \left\{ z=x + \ai y \left|\, x=\tp{(x_1,\dots,x_{p+q})}, y=\tp{(y_1,\dots,y_{p+q})} \in \R^{p+q} \right. \right\},
\]
which is equipped with a symplectic form $\omega$ given by
\begin{equation}
\label{e:omega}
  \omega(z,w)=\impart( z^* I_{p,q} w)  \qquad (z,w \in W) 
\end{equation}
Then $G$ acts on $(W,\omega)$ symplectically via $z \mapsto g z$ (matrix multiplication) for $z \in W$ and $g \in G$.
Furthermore, the action of $G$ on $(W,\omega)$ is Hamiltonian, i.e.,~ there exists a moment map $\mu: W \to \g_0^*$, whose definition we briefly recall:
%
%
if, in general, a Lie group $G$ acts on a symplectic manifold $(M,\omega)$ symplectically, a smooth $G$-equivariant map $\mu: M \to \g_0^*$ that satisfies
\begin{equation}
\label{e:definition of moment map}
 \rmd \, \< \mu, X \> =\iota(X_M) \omega \qquad \text{for all} \; X \in \g_0, 
\end{equation}
is called a moment map,
where $\iota$ stands for the contraction and $X_M$ denotes the vector field on $M$ given by
\begin{equation}
\label{e:vector field}
 X_M(p) = \left. \frac{\rmd}{\rmd t} \right|_{t=0} \exp(-t X).p    \qquad (p \in M).
\notag
\end{equation}

Under the identification that $e_i:=\tp{(0,\dots,\overset{i}1,\dots,0)} \leftrightarrow \pd_{x_i}$ and $\ai e_i \leftrightarrow \pd_{y_i}$ for $i=1,2,\dots,p+q$,
the symplectic form $\omega$ given in \eqref{e:omega} can be rewritten as
\[
  \omega = \sum_{i=1}^{p+q} \epsilon_i \rmd x_i \wedge \rmd y_i
\]
with $\epsilon_i=1$ for $i \in [p]$ and $\epsilon_{p+i}=-1$ for $i \in [q]$.

%
%
\begin{proposition}
The action of $G=\rmO(p,q)$ on $(W,\omega)$ is Hamiltonian, and the moment map $\mu:W \to \g_0^* \simeq \g_0$ is given by
\begin{equation}
\label{e:moment map}
\begin{aligned}
  \mu(z) 
    &= -\frac{\ai}{2}\left( z z^* - \tp{(z z^*)} \right) I_{p,q}
          \\
    &= (- x \tp{y} + y \tp{x}) I_{p,q}
          \\
    &= \begin{bmatrix}
        - x' \tp{y'} + y' \tp{x'}    & x' \tp{y''} - y' \tp{x''}  \\
        - x'' \tp{y'} + y'' \tp{x'}  & x'' \tp{y''} - y'' \tp{x''}
      \end{bmatrix}
\end{aligned}
\notag
\end{equation}
for $z = x + \ai y \in W$ with $x=(x',x''), \; y=(y',y'') \in \R^{p+q}$ and $x', y' \in \R^p, x'',y'' \in \R^q$.
\end{proposition}
\begin{proof}
See e.g.~\cite[{Proposition 1.4.6}]{CG97}.  
\end{proof}

\begin{remark}
Recall that the moment map $\mu_{\U}: W \to \mathfrak{u}(p,q)^* \simeq \mathfrak{u}(p,q)$ for the action of $\U(p,q)$ on $(W,\omega)$ is given by
\[
  \mu_{\U}(z) = - \ai z z^* I_{p,q}  \qquad (z = x + \ai y \in W)
\]
where we identify $\mathfrak u(p,q)^*$ with $\mathfrak u(p,q)$ via the invariant bilinear form $B$ given by $B(X,Y)=(1/2)\trace{X Y}$.
Therefore, the moment map $\mu$ in the proposition is related to $\mu_{\U}$ by 
\[
   \mu(z) = \frac{ \mu_{\U}(z) + \overline{\mu_{\U}(z)} }{2}.
\]
Namely, one has $\mu = \repart \mu_{\U}$.
\end{remark}

We define a Poisson bracket by 
\[
  \{f,g\} = \omega(\xi_g,\xi_f),
\]
where $\xi_f$ denotes the Hamiltonian vector field on $W$ corresponding to $f \in C^\infty(W)$, i.e.~the vector field that satisfies $\iota(\xi_f) \omega = \rmd f$.
Then the Poisson bracket among the coordinate functions are given by
\begin{equation}
  \{ x_i,y_j \} = - \delta_{i,j} \epsilon_i, \qquad \{ x_i,x_j \} = \{ y_i, y_j \}=0
\notag
\end{equation}  
for $i,j=1,2,\dots,p+q$.
The Dirac's quantization conditions requires that
\[
 \{f_1,f_2\} = f_3 \quad \text{implies} \quad [\hat f_1,\hat f_2] = -\ai \hbar \hat f_3 
\] 
for $f_i \in C^{\infty}(W)$ (see e.g.~\cite{Woodhouse92}). 
Thus, we quantize the coordinate functions as follows:
\begin{equation}
\label{e:CQ}
\begin{aligned}
  \hat x_i &= x_i,                       & \qquad  \hat y_i &= -\ai \hbar \pd_{x_i},   & \quad (& i=1,\dots,p), \\
  \hat x_{\jbar} &= -\ai \hbar \pd_{y_{\jbar}}, & \qquad \hat y_{\jbar} &= y_{\jbar},  & \quad (& j=1,\dots,q),
\end{aligned}
\end{equation}
where $\pd_{x_i}$ and $\pd_{y_{\jbar}}$ denote $\pd / \pd x_i$ and $\pd / \pd y_{\jbar}$ respectively.
In what follows, we set $\hbar=1$ for brevity.

The quantization \eqref{e:CQ} corresponds to a Lagrangian subspace $V$ of $W$ given by
\begin{equation}
\label{e:lagrangian}
 V = \< e_1, \dots, e_p, \ai e_{\bar 1}, \dots, \ai e_{\bar q} \>_{\R}
\end{equation}
in the sense that the quantized operators are realized in $\PD(V)$, the ring of polynomial coefficient differential operators on $V$.
Therefore, 
the quantized moment map $\hat \mu$ is given by
\begin{equation}
  \hat \mu = (-\hat x \tp{\hat y} + \hat y \tp{\hat x}) I_{p,q}
       =\begin{bmatrix}
          \ai(x' \tp{\pd_{x'}} - \pd_{x'} \tp{x'})  & x' \tp{y''} + \pd_{x'} \tp{\pd_{y''}}  \\
	  \pd_{y''} \tp{\pd_{x'}} + y'' \tp{x'}   & \ai(y'' \tp{\pd_{y''}} - \pd_{y''} \tp{y''}) 
	\end{bmatrix},
\notag
\end{equation}
where 
\begin{equation}
 \begin{aligned}
  \hat x &= \tp{(\hat x_1,\dots,\hat x_{p+q})}=\tp{(x',-\ai \pd_{y''})}, \\
  \hat y &= \tp{(\hat y_1,\dots,\hat y_{p+q})}=\tp{(-\ai \pd_{x'},y'')},  
 \end{aligned} 
\notag
\end{equation}
and
\begin{equation}
\begin{aligned}
 x' &=\tp{(x_1,\dots,x_p)},               & \qquad  \pd_{x'} &= \tp{(\pd_{x_1},\dots, \pd_{x_p})}, \\
 y''&=\tp{(y_{\bar 1},\dots,y_{\bar q})}, & \qquad  \pd_{y''} &= \tp{(\pd_{y_{\bar 1}},\dots,\pd_{y_{\bar q}})}.
\end{aligned}
\notag
\end{equation}
Note that $x_1,\dots,x_p,y_{\bar 1},\dots,y_{\bar q}$ are considered to be the coordinate functions on $V$ 
with respect to the basis $e_1,\dots,e_p, \ai e_{\bar 1},\dots, \ai e_{\bar q}$.

%
%

\begin{theorem}
For $X \in \g$, set $\pi(X):=\ai \< \hat \mu, X \>$.
Then $\pi: \g \to \PD(V)$ is a Lie algebra homomorphism.
In terms of the basis \eqref{e:basis}, it is given by  
\begin{equation}
\label{e:varpi}
  \pi(X) = \begin{cases}
     -x_j \pd_{x_i} + x_i \pd_{x_j}  &\text{if} \quad X=X_{i,j}^+ ; \\
     -y_{\jbar} \pd_{y_{\ibar}} + y_{\ibar} \pd_{y_{\jbar}}  &\text{if} \quad X=X_{\ibar,\jbar}^+ ; \\
     \ai ( x_i y_{\jbar} + \pd_{x_i} \pd_{y_{\jbar}})   &\text{if} \quad X=X_{i,j}^-.
           \end{cases}
\end{equation}
\end{theorem}
\begin{proof}
This is proved in the same manner as \cite[Theorem 2.3]{Harm_oscil_kjm} (or, one can verify the commutation relations by direct calculation).
\end{proof}


There is another canonical quantization that corresponds to the same Lagrangian subspace $V$ of $W$ as given in \eqref{e:lagrangian}.
Namely, if we quantize the coordinate functions as
\begin{equation}
\label{e:CQ2}
\begin{aligned}
  \hat x_i &= x_i,                       & \qquad \hat y_i &= -\ai \pd_{x_i},            & \quad (& i=1,\dots,p), \\
  \hat x_{\jbar} &=y_{\jbar},            & \qquad \hat y_{\jbar} &= \ai \pd_{y_{\jbar}}, & \quad (& j=1,\dots,q),
\end{aligned}
\end{equation}
then the quantized moment map, which we denote by $\hat{\mu}^{\sharp}$, is given by
\begin{equation}
  \hat{\mu}^{\sharp} = (-\hat x \tp{\hat y} + \hat y \tp{\hat x}) I_{p,q}
           = \ai 
            \begin{bmatrix}
             x' \tp{\pd_{x'}} - \pd_{x'} \tp{x'}   & x' \tp{\pd_{y''}} + \pd_{x'} \tp{y''}  \\
	     y'' \tp{\pd_{x'}} + \pd_{y''} \tp{x'}   & y'' \tp{\pd_{y''}} - \pd_{y''} \tp{y''} 
	    \end{bmatrix},
\notag
\end{equation}
where 
\begin{equation}
 \begin{aligned}
  \hat x &= \tp{(\hat x_1,\dots,\hat x_{p+q})}=\tp{(x',y'')}, \\
  \hat y &= \tp{(\hat y_1,\dots,\hat y_{p+q})}=\tp{(-\ai \pd_{x'}, \ai \pd_{y''})}.
 \end{aligned} 
\notag
\end{equation}
Hence one obtains a representation $\pi^{\sharp}: \g \to \PD(V)$ if one sets $\pi^{\sharp}(X):= \ai \<\hat \mu^{\sharp}, X\>$ for $X \in \g$. 
It is given in terms of the basis \eqref{e:basis} by
\begin{equation}
\label{e:varpi_sharp}
 \pi^{\sharp}( X ) 
  = \begin{cases}
      -x_j \pd_{x_i} + x_i \pd_{x_j}  & \text{if} \quad X=X_{i,j}^+ ; \\
      -y_{\jbar} \pd_{y_{\ibar}} + y_{\ibar} \pd_{y_{\jbar}}  & \text{if} \quad X=X_{\ibar,\jbar}^+ ; \\
      -( x_i \pd_{y_{\jbar}} + y_{\jbar} \pd_{x_i})   & \text{if} \quad X=X_{i,j}^-.
    \end{cases}
\end{equation}

%
%
\begin{remark}
\label{r:relation_with_Howe-Tan}
(i)\;
Comparing \eqref{e:CQ2} with \eqref{e:CQ}, 
one sees that $\pi^\sharp$ is related to $\pi$ through the partial Fourier transform on $\R^{p+q}$ with respect to the variables $y_{\bar 1}, \dots, y_{\bar q}$. 
In fact, if we denote the dual variable of $y_{\jbar}$ by $\eta_{\jbar}$, $j=1,2,\dots,q$, then $\pi$ and $\pi^\sharp$ interchange with each other under the correspondence
\begin{equation}
  -\ai \pd_{y_{\jbar}} \longleftrightarrow \eta_{\jbar},  \quad y_{\jbar} \longleftrightarrow \ai \pd_{\eta_{\jbar}}
     \qquad (j=1,\dots,q); 
\notag
\end{equation}
the former operators $-\ai \pd_{y_{\jbar}}$ and $\eta_{\jbar}$ are the realizations of $\hat x_{\jbar}$, 
while the latter operators $y_{\jbar}$ and $\ai \pd_{\eta_{\jbar}}$ are the realizations of $\hat y_{\jbar}$.

(ii)\;
Recall that one can obtain $\pi^\sharp$ by differentiating the left regular representation of $G=\rmO(p,q)$ on $C^{\infty}(V)$, the space of complex-valued smooth functions on $V$,
where $G$ acts on $V$ by matrix multiplication under the identification of $V$ with $\R^{p+q}$ given by $\tp{(x',\ai y'')} \leftrightarrow \tp{(x',y'')}$
(see e.g.~\cite{RS80, HT93}).
As one can see from \eqref{e:varpi} and \eqref{e:varpi_sharp}, $\pi^{\sharp}(X)$ coincides with $\pi(X)$ for all $X \in \mathfrak k$.
Thus, the action $\pi$ restricted to $\mathfrak k_0$ lifts to the action of $K$ on $C^\infty(V)$.
\end{remark}

\section{Dual Pair $( \rmO(p,q), \mathfraksl_{2}(\R) )$}
\label{s:dual_pair_in_this_paper}

Henceforth, let us denote $x'=\tp{(x_1,\dots,x_p)}$ and $y''=\tp{(y_{\bar 1},\dots,y_{\bar q})}$ by 
\begin{equation}
  x = \tp{(x_1,\dots,x_p)} \quad \text{and} \quad  y = \tp{(y_1,\dots,y_q)}
\notag
\end{equation}
respectively for the sake of simplicity if there exists no risk of confusion.
Namely, we regard $(x_1,\dots,x_p)$ and $(y_1,\dots,y_q)$ as the canonical coordinate functions on $\R^p$ and on $\R^q$ respectively.

If we denote the Casimir elements of $\g$ by $\Omega_{\g}$, 
then the corresponding Casimir operator is given by
\begin{equation}
\label{e:Casimir_opq} 
\begin{aligned}
  \pi( \Omega_{\g} ) 
     &= (\Eu_x-\Eu_y)^2 + (p-q) (\Eu_x-\Eu_y)  - 2 (\Eu_x + \Eu_y)  
          \\ 
     & \hspace{1em} - \left( r_x^2 \, r_y^2 + r_x^2 \, \Delta_x + r_y^2 \, \Delta_y + \Delta_x \, \Delta_y \right) - p q,
\end{aligned}    
\end{equation}
where 
\begin{equation}
\label{e:TDS_in_x_and_y}
\begin{aligned}
 \Eu_x &= \sum_{i \in [p]} x_i \pd_{x_i},  & \quad   r_x^2 &= \sum_{i \in [p]} x_i^2,  & \quad  \Delta_x &= \sum_{i \in [p]} \pd_{x_i}^2,  \\
 \Eu_y &= \sum_{j \in [q]} y_j \pd_{y_j},  & \quad   r_y^2 &= \sum_{j \in [q]} y_j^2,  & \quad  \Delta_y &= \sum_{j \in [q]} \pd_{y_j}^2.
\end{aligned}
\end{equation}
Now, taking account of the fact that our realization of the representation operators of $\g$ given in \eqref{e:varpi}
is a partial Fourier transform of the ones given in \cite{RS80,HT93} as we mentioned in Remark \ref{r:relation_with_Howe-Tan} (i) above,
we define elements $H,X^{+},X^{-}$ of $\PD(V)$ by
\begin{equation}
\label{e:TDS}
 H = -\Eu_x - \frac{p}2 + \Eu_y +\frac{q}2, \quad
 X^{+} = -\frac12 ( \Delta_x + r_y^2 ), \quad
 X^{-} = \frac12 ( r_x^2 + \Delta_y ). 
\end{equation}
Then, it is immediate to see that the commutation relations among them are given by
\begin{equation}
 [H, X^{+}] = 2 X^{+}, \quad [H, X^{-}] = -2 X^{-}, \quad [X^{+}, X^{-}] = H.
\notag
\end{equation}

%
%
\begin{proposition}
Let $\g':=\Cspan{ H,X^{+},X^{-} }$. 
Then $\g'$ is a Lie subalgebra of $\PD(V)^{\g}$ isomorphic to $\mathfraksl_2$ ($= \mathfraksl_2(\C)$),
where $\PD(V)^{\g}$ denotes the commutant of $\g$ in $\PD(V)$.
\end{proposition}
\begin{proof}
Note that $\pi(X_{i,j}^+)$, $i,j \in [p]$, span the Lie subalgebra isomorphic to $\mathfrak{o}_p$ commuting with $E_x, \Delta_x$ and $r_x^2$,
and that $\pi(X_{\ibar,\jbar}^+)$, $i,j \in [q]$, span the Lie subalgebra isomorphic to $\mathfrak{o}_q$ commuting with $E_y, \Delta_y$ and $r_y^2$.
Hence, it remains to show that each $\pi(X_{i,j}^{-})$ commutes with $H, X^{+}$ and $X^{-}$ given in \eqref{e:TDS}

We will only show here that $[\pi(X_{i,j}^-), X^{+}] = 0$. 
The other cases can be shown similarly.
Now, one sees
\begin{align*}
  -2 \ai [\pi(X_{i,j}^-), X^{+}]
    &= \left[ x_i y_j+\pd_{x_i} \pd_{y_j},  - \Delta_x - r_y^2 \right]
        \\
    &= \sum_{k=1}^p \left[ \pd_{x_k}^2,  x_i \right] \! y_j - \sum_{l=1}^q \pd_{x_i} \! \left[ \pd_{y_j}, y_l^2 \right] 
        \\
    &= \sum_{k=1}^p 2 \delta_{k,i} \pd_{x_k} y_j - \sum_{l=1}^q 2 \pd_{x_i} \delta_{j,l} y_l 
        \\
    &= 2 \pd_{x_i} y_j - 2 \pd_{x_i} y_j 
     = 0.
\end{align*}
This completes the proof.
\end{proof}
If one denotes the Casimir element of $\g'$ by $\Omega_{\g'}$, then the corresponding Casimir operator that is defined by
\begin{equation}
\label{e:definition_of_Casimir_sl2}
{\allowdisplaybreaks
\begin{aligned}
  \pi(\Omega_{\g'}) 
    &= H^2 + 2( X^{+} X^{-} + X^{-} X^{+}) \\
    &= H^2 -2 H + 4 X^{+} X^{-}   \\
    &= H^2 + 2 H + 4 X^{-} X^{+}
\end{aligned}
}
\notag
\end{equation}
is concretely written in terms of the operators given by \eqref{e:TDS_in_x_and_y} as follows:
\begin{equation}
\label{e:Casimir_sl2}
\begin{aligned}
  \pi(\Omega_{\g'}) 
     & = (\Eu_x-\Eu_y)^2 + (p-q) (\Eu_x-\Eu_y) - 2 (\Eu_x + \Eu_y)  \\ 
     & \hspace{1em} - \left( r_x^2 \, r_y^2 + r_x^2 \, \Delta_x + r_y^2 \, \Delta_y + \Delta_x \, \Delta_y \right) 
                        + \frac14 ( p - q )^2  - ( p + q ).
\end{aligned}
\end{equation}
It follows from \eqref{e:Casimir_opq} and \eqref{e:Casimir_sl2} that
\begin{equation}
  \pi(\Omega_{\g}) = \pi(\Omega_{\g'}) - \frac14 \left( p+q \right)^2 + \left( p+q \right) 
\notag
\end{equation}
(see \cite{Howe79, RS80}).

In what follows, we denote by $\Harm^k(\R^n)$ the space of homogeneous harmonic polynomials on $\R^n$ of degree $k$.
It is well known that $\Harm^k(\R^n)$ is an irreducible $\rmO(n)$-module and its dimension is given by

\begin{align}
  \dim \Harm^k(\R^n) 
     &= \binom{k+n-1}{n-1} - \binom{k+n-3}{n-1} 
\notag
          \\[3pt] 
     &= \frac{(k+n-3)!}{k! \, (n-2)!} \, (2k+n-2)
\notag
\end{align}
if $n \geqsl 2$ and $k \in \N$, where $\binom{\nu}{i}$ denotes the binomial coefficient.
Note that it can be further rewritten as 
\begin{equation}
\label{e:dim_Harm}
  \dim \Harm^k(\R^n) = \frac{2(k+n/2-1)}{(n-2)!} \, (k+1)(k+2) \cdots (k+n-3).
\end{equation}

Now, we will find a highest weight vector with respect to the $\g'$-action \eqref{e:TDS}, i.e.~a function $f$ on $V$ which satisfies 
\begin{equation}
\label{e:highest_weight_vec}
 H f = \lambda f
     \quad \text{and} \quad 
 X^{+} f = 0
\end{equation}  
for some $\lambda \in \C$.
Taking account of the fact that the algebra of polynomial functions on $V$, say $\PV$, can be written as
\begin{align*}
 \PV  & \simeq \C[x_1,\dots,x_p] \otimes \C[y_1,\dots,y_q] 
           \\
      & \simeq \bigoplus_{k=0}^{\infty} \left( \C [ r_x^{2} ]  \otimes \Harm^k(\R^p) \right) 
                  \otimes \bigoplus_{l=0}^{\infty} \left( \C [ r_y^{2} ] \otimes \Harm^l(\R^q) \right) 
           \\
      & \simeq  \bigoplus_{k,l=0}^{\infty}  \Harm^k(\R^p) \otimes \Harm^l(\R^q) \otimes \C [ r_x^2, r_y^2 ],
\end{align*} 
we will seek for a function that satisfies \eqref{e:highest_weight_vec} of the form 
\begin{equation}
\label{e:h1_tensor_h2_phi}
 f(x,y) = h_1(x) h_2(y) \phi(r_x^2,r_y^2), 
\end{equation}
where $h_1 \in \Harm^k(\R^p), h_2 \in \Harm^l(\R^q)$, and $\phi(s,t) \in \C[[s,t]]$
(\textbf{Caution}: we do \emph{not} assume that $\phi$ is a \emph{polynomial}). 
Namely, our function $f$ on $V$ lives in the space $\tilde \ourE$ defined by
\begin{equation}
\label{e:ambient_function_space}
  \tilde \ourE := \bigoplus_{k,l=0}^\infty {\Harm}^k(\R^p) \otimes {\Harm}^l(\R^q) \otimes \C[[r_x^2,r_y^2]]
     \qquad 
   (\text{algebraic direct sum}).
\end{equation}
Recall that the action $\pi$ of $\mathfrak k_0$ lifts to the action of $K$ on $\tilde \ourE$ as we mentioned in Remark \ref{r:relation_with_Howe-Tan} (ii),
which we denote by the same letter $\pi$.

%
%
\begin{lemma}
\label{l:laplacian_applied_to_product}
Let $\Delta=\sum_{i=1}^n \pd_{x_i}^2$ and $r^2=\sum_{i=1}^n x_i^2$.
For $h=h(x_1,\dots,x_n)$ a homogeneous harmonic polynomial on $\R^n$ of degree $d$ and for $\varphi(u)$ a smooth function in a single variable $u$,
we have
\begin{equation}
\label{e:laplacian_applied_to_product}
  \Delta ( h \varphi(r^2) ) = (4d + 2n) h \varphi'(r^2) + 4 r^2 h \varphi''(r^2).
\notag
\end{equation}
\end{lemma}
\begin{proof}
Since $\pd_{x_i} \varphi(r^2) = 2 x_i \varphi'(r^2)$ and $\pd_{x_i}^2 \varphi(r^2) = 2 \varphi'(r^2) + 4 x_i^2 \varphi''(r^2)$, one obtains
\[
  \Delta \varphi(r^2) 
    = 2 n \varphi'(r^2) + 4 r^2 \varphi''(r^2).
\]
Thus,
\begin{align*}
  \Delta ( h \varphi(r^2) ) 
     &= \sum_{i=1}^n \left( \pd_i^2 h \cdot \varphi(r^2) + 2 \pd_i h \cdot \pd_i \varphi(r^2) + h \cdot \pd_i^2 \varphi(r^2) \right) 
          \\
     &= 4 d h \varphi'(r^2) + h \Delta \varphi(r^2)
          \\
     &= 4 d h \varphi'(r^2) + h \left( 2 n \varphi'(r^2) + 4 r^2 \varphi''(r^2) \right)
          \\
     &= (4 d + 2 n) h \varphi'(r^2) + 4 r^2 h \varphi''(r^2). 
\qedhere
\end{align*}
\end{proof}

For $h_1 \in \Harm^k(\R^p)$ (resp. $h_2 \in \Harm^l(\R^q)$) given, we define its shifted degree by $\kappa_{+}(h_1):=k+p/2$ (resp. $\kappa_{-}(h_2):=l+q/2$),
which we denote just by $\kappa_{+}$ (resp. $\kappa_{-})$ if there is no risk of confusion.

It follows from Lemma \ref{l:laplacian_applied_to_product} that if $f$ is of the form in \eqref{e:h1_tensor_h2_phi} then
\begin{align*}
 X^{+} f 
   &= -\frac12 ( \Delta_x (h_1 h_2 \phi) + r_y^2 h_1 h_2 \phi ) \\
   &= -2 h_1 h_2 \left( r_x^2 \, (\pd_s^2 \phi)(r_x^2,r_y^2) + \kappa_{+} (\pd_s \phi)(r_x^2,r_y^2) + r_y^2 \, \phi(r_x^2,r_y^2) \right),
\end{align*}
which shows that $f=h_1(x) h_2(y) \phi(r_x^2,r_y^2)$ satisfies $X^{+} f=0$ if and only if $\phi$ is a solution to a differential equation
\begin{equation}
\label{e:DE_for_extremal_vec}
  s \pd_s^2 \phi + \kappa_{+} \pd_s \phi + t \phi = 0
\end{equation}
with $\kappa_{+} = \kappa_{+}(h_1) = k + p/2$. 
Solving the differential equation \eqref{e:DE_for_extremal_vec} by power series, one obtains that
\begin{equation}
\label{e:sol4DE_extremal_vec}
 \phi(s,t) = a_0 \sum_{n=0}^{\infty} \frac{(-1)^n}{n! (\kappa_{+})_n}\left( \frac{s \, t}{4} \right)^n,  
\end{equation}
where $a_0$ is an arbitrary formal power series in $t$.
Note that if one defines a power series $\Psi_\alpha$ by
\begin{equation}
\label{e:normalized_Bessel}
 \Psi_\alpha (u) 
     := \sum_{n=0}^\infty \frac{(-1)^n}{n! \, (\alpha)_n} u^n
     = 1 - \frac{u}{\alpha} + \frac{u^2}{2! \, \alpha(\alpha+1)} - \frac{u^3}{3! \, \alpha(\alpha+1)(\alpha+2)} + \cdots    
\end{equation}
for $\alpha \in \C \setminus (- \N )$,
then it converges on the whole $\C$ and is a unique solution to a differential equation
\begin{equation}
\label{e:DE_for_normalized_Bessel}
  u \Psi_\alpha''(u) + \alpha \Psi_\alpha'(u) + \Psi_\alpha(u) = 0
\end{equation}
that satisfies the initial condition $\Psi_\alpha(0)=1$.

In the sequel, we set
\begin{equation}
\label{e:psi_alpha}
 \psi_\alpha^{(n)} := \Psi_\alpha^{(n)} ( {r_x^2 r_y^2}/{4} )  
     \qquad (n \in \N)
\end{equation}
for brevity, where $\Psi_\alpha^{(n)}(u)$ denotes the $n$-th derivative of $\Psi_\alpha(u)$ in $u$.

If, in addition, $f$ satisfies that $H f = \lambda f$ for some $\lambda \in \C$, 
then the factor $a_0$ in \eqref{e:sol4DE_extremal_vec} is equal to $t^{\mu_{-}}$ up to a constant multiple with 
\( 
  \mu_{-} = (1/2) (\lambda + \kappa_{+} - \kappa_{-}) \in \N,
\)
$\kappa_{+} = \kappa_{+}(h_1)$ and $\kappa_{-} = \kappa_{-}(h_2)$.

Thus, for $h_1 \in \Harm^k(\R^p)$ and $h_2 \in \Harm^l(\R^q)$ given, a highest weight vector $f=f(x,y)$ with respect to the $\g'$-action,
i.e.~ a function that satisfies \eqref{e:highest_weight_vec} is given by
\begin{equation}
\label{e:h_weight_vec}
  f = h_1(x) h_2(y) r_y^{2 \mu_{-}} \psi_{\kappa_{+}}
\notag
\end{equation}
with
\begin{equation}
\label{e:h_weight}
  \lambda = -\kappa_{+} + \kappa_{-} + 2 \mu_{-}
     \quad (\mu_{-} \in \N).
\end{equation}

Similarly, for $h_1 \in \Harm^k(\R^p)$ and $h_2 \in \Harm^l(\R^q)$ given, 
one can show that a lowest weight vector $f=f(x,y)$ of the form \eqref{e:h1_tensor_h2_phi} with respect to the $\g'$-action, i.e.~ a function that satisfies
\begin{equation}
\label{e:lowest_weight_vec}
 H f = \lambda f
  \quad \text{and} \quad
 X^{-} f = 0
\end{equation}  
for some $\lambda \in \C$, is given by
\begin{equation}
\label{e:l_weight_vec}
  f(x,y) = h_1(x) h_2(y) r_x^{2 \mu_{+}} \psi_{\kappa_{-}}  
\notag
\end{equation}
with
\begin{equation}
\label{e:l_weight}
 \lambda = -\kappa_{+} + \kappa_{-} - 2 \mu_{+}
    \quad (\mu_{+} \in \N).
\end{equation}

%
%
Let us summarize the above argument in the following.
\begin{proposition}
\label{p:extremal_vectors}
Given $h_1 \in \Harm^k(\R^p)$ and $h_2 \in \Harm^l(\R^q)$, let $f=f(x,y)$ be a function of the form given in \eqref{e:h1_tensor_h2_phi}.
\begin{list}{}{\setlength{\itemindent}{0pt} \setlength{\leftmargin}{6pt}}
  \item[$(1)$]
If $f$ is a highest weight vector satisfying \eqref{e:highest_weight_vec} with respect to the $\g'$-action, then it is given by
\begin{equation}
  f(x,y) = h_1(x) h_2(y) r_y^{2 \mu_{-}} \psi_{\kappa_+}
\notag
\end{equation}
with $\lambda = -\kappa_{+} + \kappa_{-} + 2 \mu_{-}$.
Moreover, such a function is unique up to a constant multiple.

  \item[$(2)$]
If $f$ is a lowest weight vector satisfying \eqref{e:lowest_weight_vec} with respect to the $\g'$-action, then it is given by
\begin{equation}
  f(x,y) = h_1(x) h_2(y) r_x^{2 \mu_{+}} \psi_{\kappa_-}
\notag
\end{equation}
with $\lambda = -\kappa_{+} + \kappa_{-} - 2 \mu_{+}$.
Moreover, such a function is unique up to a constant multiple.
\end{list}
Here $\kappa_+ = \kappa_{+}(h_1) = k + p/2$, $\kappa_- = \kappa_{-}(h_2) = l + q/2$, $\mu_+, \mu_- \in \N$, 
and $\psi_{\kappa_{\pm}}$ is an element of $\C[[ r_x^2, r_y^2]]$ given by \eqref{e:psi_alpha} with $\alpha=\kappa_{\pm}$ and $n=0$.
\end{proposition}

Taking account of the discussion so far, let us introduce the subspace $\ourE$ of $\tilde \ourE$ by 
\begin{equation}
\label{e:ourE}
  \ourE := \Cspan{ h_1(x) h_2(y) \rho_x^a \rho_y^b \psi_\alpha \in \tilde \ourE
            \left|  
             \begin{array}[c]{l}
	      h_1 \in \Harm^k(\R^p), h_2 \in \Harm(\R^q), \\
              a,b \in \N, \alpha \in \C \setminus (-\N)
	     \end{array} 
            \right.
                }.
\notag
\end{equation}
Then one will find that $\ourE$ is stable under the action of $(\g,K)$ as well as that of $\g'$ (see Propositions \ref{p:sl2-action_on_ourE} and \ref{p:p-action} below).

%
%
\begin{remark}
(i)\;
The function $\Psi_\alpha$ given in \eqref{e:normalized_Bessel} can be written in terms of the Bessel function $J_\nu$ of the first kind of order $\nu$
\[
  J_\nu (t)  = \sum_{n=0}^{\infty} \frac{(-1)^n}{n! \, \Gamma(n+ \nu + 1)} \left( \frac{t}{2} \right)^{\nu + 2 n}
\]
that solves the Bessel's differential equation
\begin{equation}
\label{e:DE_Bessel}
  \frac{\rmd^2 w}{\rmd \, t^2} + \frac{1}{t} \frac{\rmd w}{\rmd \, t} + \left( 1 - \frac{\nu^2}{t^2} \right) w = 0
\end{equation}
(see e.g.~\cite{WW27}).
Namely, one has
\begin{equation}
\label{e:Psi_and_Bessel_fun}
  \Psi_{\alpha}(u) = \Gamma( \alpha ) \, u^{-(\alpha-1) / 2} J_{\alpha-1}( 2 \, {u}^{1/2} ).
\end{equation}
Therefore,
\begin{equation}
 \psi_\alpha = \Gamma(\alpha) \left( \frac{r_x r_y}{2} \right)^{-(\alpha - 1)} J_{\alpha-1} (r_x r_y).
\notag
\end{equation}
Note that \eqref{e:DE_for_normalized_Bessel} corresponds to \eqref{e:DE_Bessel} under \eqref{e:Psi_and_Bessel_fun}.

%
%
(ii)\;
Recall that our representation $\pi$ is related to $\pi^\sharp$ via the partial Fourier transform with respect to $y_1,\dots,y_q$, 
as we mentioned in Remark \ref{r:relation_with_Howe-Tan} (i).
Namely, one can obtain $\pi^\sharp$ by replacing $-\ai \pd_{y_j}$ and $y_j$ in $\pi$ by $\eta_j$ and $\ai \pd_{\eta_j}$, $j=1,\dots,q$, respectively.
Under this correspondence, one finds that $H = -\Eu_x - p/2 + \Eu_y + q/2$ and $X^{+} = -\frac12( \Delta_x + r_y^2)$ transforms, up to constant multiples,
into the shifted degree operator $\tilde\Eu_{p,q}$ and the d'Alembertian $\Box_{p,q}$ on $\R^{p+q}$ that are given by
\begin{align*}
  \tilde\Eu_{p,q} &= \sum_{i=1}^p x_i \pd_{x_i} + \sum_{j=1}^q \eta_j \pd_{\eta_j} + \frac{p-q}{2},
          \\ 
  \Box_{p,q} &= \sum_{i=1}^{p} \pd_{x_i}^2 - \sum_{j=1}^{q} \pd_{\eta_j}^2,
\end{align*}
respectively.
Therefore, the highest weight vector $f$ satisfying $H f = \lambda f$ for some $\lambda \in \C$ and $X^{+} f = 0$ 
corresponds to a homogeneous solution $\tilde f$ to the equation $\square_{p,q} \tilde f=0$.
\end{remark}

Note that $\Psi_{\alpha}^{(n)}$ is equal to $\Psi_{\alpha+n}$ up to a constant multiple.
In fact, differentiating both sides of \eqref{e:DE_for_normalized_Bessel} $n$ times, one obtains
\begin{equation}
\label{e:DE_for_n-th_derivative_Bessel_fnc}  
  u \Psi_\alpha^{(n+2)} (u) + (\alpha + n) \Psi_\alpha^{(n+1)} (u) + \Psi_\alpha^{(n)} (u) = 0.
\end{equation}
Since $\Psi_{\alpha+n}$ is a unique solution to \eqref{e:DE_for_normalized_Bessel} with $\alpha$ replaced by $\alpha+n$ that satisfies $\Psi_{\alpha+n}(0)=1$,
it follows that $\Psi_\alpha^{(n)} = {(-1)^n}/{(\alpha)_n} \Psi_{\alpha+n}$.
Thus, one obtains
\begin{equation}
\label{e:n-th_derivative_psi}
 \psi_\alpha^{(n)} = \frac{(-1)^n}{(\alpha)_n} \psi_{\alpha+n}
     \qquad  (n \in \N).
\end{equation}
In what follows, we set $\rho_x:=r_x^2/2$ and $\rho_y:=r_y^2/2$ for economy of space.
Then, it follows from \eqref{e:DE_for_n-th_derivative_Bessel_fnc} and \eqref{e:n-th_derivative_psi} that
\begin{equation}
\label{e:recursive_formula_psi}
  \rho_x \rho_y \psi_{\alpha+2} = \alpha (\alpha + 1) ( \psi_{\alpha+1} - \psi_\alpha )
\end{equation}
for $\alpha \in \C \setminus (-\N)$.
Furthermore, setting $\rho:=r^2/2$, one can rewrite \eqref{l:laplacian_applied_to_product} as
\begin{equation}
\label{e:laplacian_applied_to_product2}
  \frac12 \Delta ( h \varphi(\rho) ) = \left( d + \frac{n}{2} \right) h \varphi'(\rho) + h \rho \varphi''(\rho),
\end{equation}
where $h, \Delta, r^2$ and $\varphi$ are as in Lemma \ref{l:laplacian_applied_to_product}.

%
%
\begin{proposition}
\label{p:sl2-action_on_ourE}
For $f= h_1 h_2 \rho_x^a \rho_y^b \psi_{\alpha} \in \ourE$, one has
\begin{align}
  H \bigl( h_1 h_2 \rho_x^a \rho_y^b \psi_{\alpha} \bigr)
   & = ( -\kappa_{+} + \kappa_{-} -2 a + 2 b ) h_1 h_2 \rho_x^a \rho_y^b  \psi_{\alpha},
     \label{e:sl2-action_H}
          \\
  X^{+} \bigl( h_1 h_2  \rho_x^a \rho_y^b \psi_{\alpha} \bigr)
   & = h_1 h_2 \Bigl( \mbox{\small $-a ( \kappa_{+} + a - 1 )$} \rho_x^{a-1} \rho_y^b \psi_{\alpha}
        + \tfrac{\kappa_{+} + 2 a -\alpha}{\alpha} \rho_x^a \rho_y^{b+1} \psi_{\alpha+1} \Bigr),
     \label{e:sl2-action_X}
          \\
  X^{-} \bigl( h_1 h_2 \rho_x^a \rho_y^b \psi_{\alpha} \bigr)
   & = h_1 h_2 \Bigl( \mbox{\small $b ( \kappa_{-} + b - 1 )$} \rho_x^a \rho_y^{b-1}  \psi_{\alpha}
        - \tfrac{\kappa_{-} + 2 b -\alpha}{\alpha} \rho_x^{a+1} \rho_y^b \psi_{\alpha+1} \Bigr).
     \label{e:sl2-action_Y}
\end{align}
In particular, the $\g'$-action preserves the $K$-type of each element of $\ourE$.
\end{proposition}
\begin{proof}
It is immediate to show \eqref{e:sl2-action_H}, and we will only show \eqref{e:sl2-action_X} here; the other case \eqref{e:sl2-action_Y} can be shown similarly.

Setting $\varphi(u):=u^a \Psi_\alpha(\rho_y u)$, one sees
\begin{align*}
 \varphi'(u) &= a u^{a-1} \Psi_\alpha(\rho_y u) + u^a \rho_y \Psi_\alpha'(\rho_y u),
     \\
 \varphi''(u) &= a(a-1)u^{a-2} \Psi_\alpha(\rho_y u) + 2 a u^{a-1} \rho_y \Psi'_\alpha(\rho_y u) + u^a \rho_y^2 \Psi''_\alpha(\rho_y u).
\end{align*}
Hence it follows from \eqref{e:laplacian_applied_to_product2} that
\begin{align*}
 \frac12 \Delta_x (h_1 \rho_x^a \psi_\alpha)
   &= a (\kappa_{+} + a-1) h_1 \rho_x^{a-1} \psi_\alpha + ( \kappa_{+} + 2 a ) h_1 \rho_x^a \rho_y \psi'_\alpha + h_1 \rho_x^{a+1} \rho_y^2 \psi''_\alpha
     \\
   &= h_1 \left( a(\kappa_{+} + a-1) \rho_x^{a-1} \psi_\alpha + (\kappa_{+}+2 a-\alpha) \rho_x^a \rho_y \psi'_\alpha - \rho_x^a \rho_y \psi_\alpha \right)  
\end{align*}
since $\rho_x \rho_y \psi''_\alpha = - \alpha \psi'_\alpha - \psi_\alpha$.
Therefore, one obtains that
\begin{align*}
  X^+ (h_1 h_2 \rho_x^a \rho_y^b \psi_\alpha)
   &= -\frac12( \Delta_x + 2 \rho_x ) (h_1 h_2 \rho_x^a \rho_y^b \psi_\alpha)
        \\
   &= -a(\kappa_{+} + a - 1) h_1 h_2 \rho_x^{a-1} \rho_y^b \psi_\alpha - (\kappa_{+} + 2 a -\alpha) h_1 h_2 \rho_x^a \rho_y^{b+1} \psi'_\alpha,  
\end{align*}
which, by \eqref{e:n-th_derivative_psi}, equals the right-hand side of \eqref{e:sl2-action_X}.
This completes the proof.
\end{proof}

We conclude this section by calculating the action of $\mathfrak p$ on $\ourE$, i.e.~ $\pi(X_{i,j}^-) f$ for $X_{i,j}^{-} \in \mathfrak{p}$ and $f \in \ourE$.

For a homogeneous polynomial $P=P(x_1,\dots,x_n)$ on $\R^n$ of degree $d$, set
\begin{equation}
\label{e:daggered_polynom}
  P^{\dag} := P - \frac{r^2}{4 (d + n/2 - 2)} \Delta P,
\notag
\end{equation}
where $\Delta=\sum_{i=1}^n \pd_{x_i}^2$ and $r^2=\sum_{i=1}^n x_i^2$.
Note that if $\Delta^2 P=0$ then $P^\dag$ is harmonic by Lemma \ref{l:laplacian_applied_to_product},
and that if $h=h(x_1,\dots, x_n)$ is harmonic then $\Delta(x_i h)=2 \pd_{x_i} h$ and $\Delta^2(x_i h)=0$.

%
%
\begin{proposition}
\label{p:p-action}
For $f= h_1 h_2 \rho_x^a \rho_y^b \psi_{\alpha} \in \ourE$, one has
\begin{equation}
\label{e:p-action}
{\allowdisplaybreaks
\begin{aligned}
 -\ai & \pi (X_{i,j}^-) \bigl( h_1 h_2 \rho_x^a \rho_y^b \psi_\alpha \bigr)  \\
   = & \, 
     (\pd_{x_i} h_1) (\pd_{y_j} h_2) \rho_x^a \rho_y^b
        \Bigl( 
            \tfrac{ (\kappa_{+} + a - \alpha) (\kappa_{-} + b - \alpha)}{(\kappa_{+}-1)(\kappa_{-}-1)} \, \psi_{\alpha} 
          + \tfrac{(\alpha-1)(\kappa_{+} + \kappa_{-} + a + b - \alpha - 1)}{(\kappa_{+}-1)(\kappa_{-}-1)} \, \psi_{\alpha-1} 
        \Bigr)
      \\
        & 
     + (\pd_{x_i} h_1) (y_j h_2)^\dag 
        \Bigl(
          - \tfrac{\kappa_{+}+a+b-\alpha}{\alpha(\kappa_{+}-1)} \rho_x^{a+1} \rho_y^b \psi_{\alpha+1} 
          + \tfrac{b (\kappa_{+} + a - 1)}{\kappa_{+}-1} \rho_x^a \rho_y^{b-1} \psi_{\alpha}
        \Bigr) 
      \\
        & 
     + (x_i h_1)^\dag (\pd_{y_j} h_2)
        \Bigl(
          - \tfrac{\kappa_{-}+a+b-\alpha}{\alpha(\kappa_{-}-1)} \rho_x^a \rho_y^{b+1} \psi_{\alpha+1}
          + \tfrac{a (\kappa_{-} + b - 1)}{\kappa_{-}-1} \rho_x^{a-1} \rho_y^b \psi_{\alpha}
        \Bigr) 
      \\
        & 
     + (x_i h_1)^\dag (y_j h_2)^\dag
        \Bigl(
          - \tfrac{a+b+1-\alpha}{\alpha} \rho_x^a \rho_y^b \psi_{\alpha+1}
          + a b \rho_x^{a-1} \rho_y^{b-1} \psi_{\alpha} 
        \Bigr). 
\end{aligned}
}
\end{equation}
\end{proposition}
\begin{proof}
Since $\pd_{x_i} \psi_\alpha = \rho_y x_i \psi'_\alpha$ and $\pd_{y_j} \psi_\alpha = \rho_x y_j \psi'_\alpha$, one obtains
\begin{equation}
\label{e:varpi_X_{ij}^{-} f}
\begin{aligned}
 - & \ai \pi(X_{i,j}^-) f = (\pd_{x_i} \pd_{y_j} + x_i y_j) (h_1 h_2 \rho_x^a \rho_y^b \psi_\alpha) \\
   &= (\pd_{x_i} h_1) (\pd_{y_j} h_2) \rho_x^a \rho_y^b \psi_\alpha  
     + (x_i h_1) (\pd_{y_j} h_2) \left( a \rho_x^{a-1} \rho_y^b \psi_\alpha + \rho_x^a \rho_y^{b+1} \psi'_\alpha \right)  \\
   & \; + (\pd_{x_i} h_1) (y_j h_2) \left( b \rho_x^a \rho_y^{b-1} \psi_\alpha + \rho_x^{a+1} \rho_y^b \psi'_\alpha \right)  \\
   & \; + (x_i h_1) (y_j h_2) 
           \left( \rho_x^a \rho_y^b \psi_\alpha + a b \rho_x^{a-1} \rho_y^{b-1} \psi_\alpha + (a+b+1) \rho_x^a \rho_y^b \psi'_\alpha + \rho_x^{a+1} \rho_y^{b+1} \psi''_\alpha \right).
\end{aligned}  
\end{equation}
Now, by definition, one has
\begin{equation}
\label{e:x_i_h_1_y_j_h_2}
  x_i h_1 = (x_i h_1)^\dag + \frac{\rho_x}{\kappa_+ - 1} \pd_{x_i} h_1 
     \quad \text{and} \quad  
  y_j h_2 = (y_j h_2)^\dag + \frac{\rho_y}{\kappa_- - 1} \pd_{y_j} h_2.
\end{equation}
Substituting \eqref{e:x_i_h_1_y_j_h_2} into \eqref{e:varpi_X_{ij}^{-} f}, and using the relation \eqref{e:n-th_derivative_psi} and \eqref{e:recursive_formula_psi},
one sees that the coefficient of $(\pd_{x_i} h_1) (\pd_{y_j} h_2)$ in \eqref{e:varpi_X_{ij}^{-} f} equals 
the one of $(\pd_{x_i} h_1)(\pd_{y_j} h_2)$ in the right-hand side of $\eqref{e:p-action}$.
One can verify that each coefficients of $(\pd_{x_i} h_1) (y_j h_2)^\dag, (x_i h_1)^\dag (\pd_{y_j} h_2)$ and $(x_i h_1)^\dag (y_j h_2)^\dag$ in \eqref{e:varpi_X_{ij}^{-} f}
equals the one of the corresponding terms in \eqref{e:p-action} similarly.
This completes the proof.
\end{proof}

\section{$(\g,K)$-module associated with finite-dimensional $\mathfraksl_2$-module}
\label{s:(g,K)-module}

If $f \in \ourE$ satisfies $H f=\lambda f,\; X^{+} f = 0$ and $(X^{-})^{m+1} f=0$ (resp. $H f=\lambda f,\; X^{-} f=0$ and $(X^{+})^{m+1} f=0$) for some $m \in \N$, 
then it follows from the representation theory of $\g' = \mathfraksl_2$ that $\lambda = m$ (resp. $\lambda = -m$).
Thus, we introduce $(\g,K)$-modules associated with the finite-dimensional $\mathfraksl_2$-module $F_m$ of dimension $m+1$ as follows, which are the main objects of this paper.

%
%
\begin{definition}
\label{d:Module^{pm}_{m}}
Given $m \in \N$, we define $(\g,K)$-modules $\Module{\pm}{m}$ by
\begin{align}
 \Module{+}{m} 
   &:= \left\{ f \in \ourE 
                \left| 
                  \begin{array}{l}
                     H f = m f, \quad X^{+} f = 0, \\
                     (X^{-})^j f \ne 0 \; (1 \leqsl j \leqsl m), \; (X^{-})^{m+1} f=0
	          \end{array}
                \right. 
       \right\},
\notag
     \\[3pt]
 \Module{-}{m} 
   &:= \left\{ f \in \ourE
                 \left|
                   \begin{array}{l}
                     H f = -m f, \quad X^{-} f = 0, \\
                     (X^{+})^j f \ne 0 \; (1 \leqsl j \leqsl m), \; (X^{+})^{m+1} f=0
	           \end{array} 
                 \right.
       \right\}.
\notag
\end{align}
The modules $\Module{\pm}{m}$ were originally introduced in \cite{RS80} without the condition of finite dimensionality. 
Note that $\Module{+}{0}$ is identical to $\Module{-}{0}$ by definition and that both $\Module{\pm}{m}$ should correspond to the $\mathfraksl_2$-module $F_m$ under the Howe duality (cf.~\cite{Howe_remarks}).
\end{definition}

If $\Module{+}{m} \ne \{ 0 \}$ (resp. $\Module{-}{m} \ne \{ 0 \}$), then one sees that $p \equiv q \hspace{-4.5pt}\mod 2$;
for, if one takes a non-zero $f=h_1 h_2 \rho_y^{\mu_-} \psi_{\kappa_{+}} \in \Module{+}{m}$ (resp. $f=h_1 h_2 \rho_x^{\mu_+} \psi_{\kappa_-} \in \Module{-}{m}$)
with $h_1 \in \Harm^k(\R^p), h_2 \in \Harm^l(\R^q)$ and $\mu_{\pm} \in \N$, 
then
\[
  \pm m = -\kappa_{+} + \kappa_{-} \pm 2 \mu_{\mp} = -k +l -\frac{p-q}2 \pm 2 \mu_{\mp} \in \Z
\]
by \eqref{e:h_weight} (resp. \eqref{e:l_weight}).
Hence one obtains ${(p-q)}/{2} \in \Z$.
Therefore, we assume that $p \equiv q \hspace{-4.5pt}\mod 2$ in the rest of this paper.

%
%
\begin{lemma}
\label{l:HWV_sl2_and_LWV_sl2}
For $h_1 \in \Harm^k(\R^p)$, $h_2 \in \Harm^l(\R^q)$ and $m \in \N$,
let
\[
  v^+ = h_1 h_2 \rho_y^{\mu_-} \psi_{\kappa_+} \in \Module{+}{m} 
    \quad \text{and} \quad 
  v^- = h_1 h_2 \rho_x^{\mu_+} \psi_{\kappa_-} \in \Module{-}{m},
\]
where $\mu_+, \mu_- \in \N$ such that $\mu_{+} + \mu_{-} = m$.
Then the $\g'$-module generated by $v^+$ coincides with the one generated by $v^-$:
\begin{equation}
  \< v^+ \>_{\g'} = \< v^- \>_{\g'}.
 \notag
\end{equation} 
\end{lemma}
\begin{proof}
Both $v^+$ and $(X^+)^m v^-$ (resp. $v^-$ and $(X^-)^m v^+$) are elements of $\ourE \subset \tilde\ourE$ that are highest (resp. lowest) weight vectors of weight $m$ (resp. $-m$) under $\g'$-action.
Namely, they are solutions in $\tilde\ourE$ to the differential equation
\[
  H f = \pm m f \quad \text{and} \quad X^{\pm} f =0.
\]
As we mentioned Proposition \ref{p:extremal_vectors}, they are respectively equal to each other up to a constant multiple.
This completes the proof.
\end{proof}

%
%
\begin{proposition}
\label{thm_item:M_m^+_and_M_m^-_are_isom}
For $m \in \N$, $\Module{+}{m}$ and $\Module{-}{m}$ are isomorphic to each other.
\end{proposition}
\begin{proof}
For $h_1 \in \Harm^k(\R^p)$ and $h_2 \in \Harm^l(\R^q)$, set 
\[
 v^+ = h_1 h_2 \rho_y^{\mu_-} \psi_{\kappa_+} \in \Module{+}{m}
   \quad \text{and} \quad 
 v^- = h_1 h_2 \rho_x^{\mu_+} \psi_{\kappa_-} \in \Module{-}{m}. 
\]
If $\mu_+ + \mu_- = m$, then $(X^+)^m v^-$ is equal to $v^+$ up to a constant multiple as in Lemma \ref{l:HWV_sl2_and_LWV_sl2}, 
and thus, $(X^+)^m (X^-)^m v^+$ is equal to $v^+$ up to a constant multiple.
Moreover, this constant is non-zero and is independent of the $K$-type of $v^+$.
In fact, using the relations \eqref{e:TDS}, one has
\begin{align*}
  H (X^-)^i 
    &= (X^-)^i H + [H,(X^-)^i] 
          \\
    &= (X^-)^i H + [H, X^-] (X^-)^{i-1} + X^- [H, X^-] (X^-)^{i-2} + \cdots + (X^-)^{i-1} [H,X^-]
          \\
    &= (X^-)^i H -2i (X^-)^i
\end{align*}
for $i \in \N$, hence,
\begin{align*}
  (X^+)^j & (X^-)^j v^+
        \\ 
     &= (X^+)^{j-1} \left( (X^-)^j X^+ + [ X^+, (X^-)^{j} ] \right) v^+
        \\
     &= (X^+)^{j-1} \left( H (X^-)^{j-1} + X^- H (X^-)^{j-2} + \cdots + (X^-)^{j-2} H X^- + (X^-)^{j-1} H \right) v^+
        \\
     &= j (m-j+1) \, (X^+)^{j-1} (X^-)^{j-1} v^+
        \\
     &= j(m-j+1) \cdot (j-1)(m-j+2) \, (X^+)^{j-2} (X^-)^{j-2} v^+
        \\
     &= \cdots
        \\
     &= j! \ffac{m}{j} \, v^+
\end{align*}
for $j \in \N$ since $X^+ v^+ = 0$.
In particular, $(X^+)^m (X^-)^m v^+ = (m!)^2 v^+$.
Therefore, 
\[
  (X^-)^m : \Module{+}{m} \longrightarrow \Module{-}{m}
\]
provides an isomorphism of $(\g,K)$-module.
This completes the proof.
\end{proof}

Now we prepare two lemmas to prove our main result. 
Note that Lemma \ref{l:p-action_on_extremal_vectors} below is just a special case of Proposition \ref{p:p-action}.
However, we state it separately to highlight the r\^{o}le of the relation \eqref{e:recursive_formula_psi}.

%
%
\begin{lemma}
\label{l:p-action_on_extremal_vectors}
Let $h_1 \in \Harm^k(\R^p)$ and $h_2 \in \Harm^l(\R^q)$, and set $\kappa_+ = k+p/2$, $\kappa_- = l+q/2$.
\begin{list}{}{\setlength{\itemindent}{0pt} \setlength{\leftmargin}{6pt}}
  \item[$(1)$]
For a highest weight vector $f=h_1 h_2 \rho_y^{\mu_{-}} \psi_{\kappa_+} \in \ourE$, $\pi(X_{i,j}^{-}) f$ is given by
\begin{equation}
\label{e:p-action_on_HWV}
\begin{aligned}
  - \ai \pi (X_{i,j}^-) \bigl( h_1 h_2 \rho_y^{\mu_{-}} \psi_{\kappa_{+}} \bigr) 
    = & \, \frac{\kappa_{-} + \mu_{-} - 1}{\kappa_{-}-1} (\pd_{x_i} h_1) (\pd_{y_j} h_2) \rho_y^{\mu_{-}} \psi_{\kappa_{+}-1}     \\
      & +  \mu_{-} (\pd_{x_i} h_1) (y_j h_2)^\dagger \rho_y^{\mu_{-}-1} \psi_{\kappa_{+}-1}     \\
      & +  \frac{\kappa_{+}-\kappa_{-}-\mu_{+}}{\kappa_{+}(\kappa_{-}-1)} (x_i h_1)^\dag (\pd_{y_j} h_2) \rho_y^{\mu_{-}+1} \psi_{\kappa_{+}+1}   \\
      & +  \frac{\kappa_{+}-\mu_{+}-1}{\kappa_{+}} (x_i h_1)^\dagger (y_j h_2)^\dagger \rho_y^{\mu_{-}} \psi_{\kappa_{+}+1}.
\end{aligned}
\end{equation}

  \item[$(2)$]
For a lowest weight vector $f=h_1 h_2 \rho_x^{\mu_{+}} \psi_{\kappa_-}$, $\pi(X_{i,j}^{-}) f$ is given by
\begin{equation}
\label{e:p-action_on_LWV}
\begin{aligned}
   - \ai \pi (X_{i,j}^-) \bigl( h_1 h_2 \rho_x^{\mu_+} \psi_{\kappa_{-}} \bigr) 
     = & \, \frac{\kappa_{+} + \mu_{+} - 1}{\kappa_{+}-1} (\pd_{x_i} h_1) (\pd_{y_j} h_2) \rho_x^{\mu_{+}} \psi_{\kappa_{-}-1}   \\
       & + \frac{\kappa_{-}-\kappa_{+}-\mu_{+}}{\kappa_{-}(\kappa_{+}-1)} (\pd_{x_i} h_1) (y_j h_2)^\dagger \rho_x^{\mu_{+}+1} \psi_{\kappa_{-}+1}   \\
       & + \mu_{+} (x_i h_1)^\dag (\pd_{y_j} h_2) \rho_x^{\mu_{+}-1} \psi_{\kappa_{-}-1}    \\
       & + \frac{\kappa_{-}-\mu_{+}-1}{\kappa_{-}} (x_i h_1)^\dagger (y_j h_2)^\dagger \rho_x^{\mu_{+}} \psi_{\kappa_{-}+1}.
\end{aligned}
\end{equation}

\end{list}
  
\end{lemma}
\begin{proof}
We only show \eqref{e:p-action_on_HWV} here.
The other formula \eqref{e:p-action_on_LWV} can be shown similarly.

Set $a=0, b=\mu_-$ and $\alpha=\kappa_+$ in \eqref{e:p-action}.
Then, using the relation \eqref{e:recursive_formula_psi} in this case, i.e.
\[
  \rho_x \rho_y \psi_{\kappa_{+}+1} = \kappa_{+}(\kappa_{+}-1)(\psi_{\kappa_{+}} - \psi_{\kappa_{+}-1}),
\]
one sees that the coefficient of $(\pd_{x_i} h_1) (y_j h_2)^\dagger$ equals
\begin{align*}
  & - \frac{\mu_{-}}{\kappa_{+}(\kappa_{+}-1)} \rho_x \rho_y^{\mu_{-}} \psi_{\kappa_{+} + 1}
    + \mu_{-} \rho_y^{\mu_{-} -1} \psi_{\kappa_{+}}
           \\
 =& \, \mu_{-} \rho_y^{\mu_{-} - 1} \psi_{\kappa_{+}-1}.
\end{align*}  
To show for the other coefficients is trivial and omitted.
\end{proof}

%
%
\begin{lemma}
\label{l:iteration_sl2-action}
Let $h_1 \in \Harm^k(\R^p)$ and $h_2 \in \Harm^l(\R^q)$, and set $\kappa_+ = k+p/2$, $\kappa_- = l+q/2$.
\begin{list}%
  {\upshape (\arabic{mylist})}{\setlength{\itemindent}{3pt} \setlength{\leftmargin}{3pt} \usecounter{mylist}}
    \item
     For a highest weight vector $f = h_1 h_2 \rho_y^{\mu_-} \psi_{\kappa_+}$ of weight $\lambda = -\kappa_{+} + \kappa_{-} + 2 \mu_{-}$, one has
     \begin{equation}
\label{e:iteration_X^{-}}
      (X^{-})^{\nu} f = h_1 h_2 \sum_{i=0}^{\nu} \binom{\nu}{i} 
           \frac{ \ffac{-\lambda + \nu -1}{i} \ffac{\mu_{-}}{\nu - i} \ffac{\kappa_{-} + \mu_{-} -1}{\nu - i} }{\rfac{\kappa_{+}}{i}} \,
           \rho_x^i \rho_y^{\mu_{-} - \nu +i} \psi_{\kappa_{+} + i}
     \end{equation}
     for $\nu=0,1,2,\dots$.

     \item
     For a lowest weight vector $f = h_1 h_2 \rho_x^{\mu_+} \psi_{\kappa_-}$ of weight $\lambda = -\kappa_{+} + \kappa_{-} - 2 \mu_{+}$, one has
     \begin{equation}
\label{e:iteration_X^{+}}
      (X^{+})^{\nu} f = (-1)^{\nu} h_1 h_2 \sum_{i=0}^{\nu} \binom{\nu}{i} 
              \frac{ \ffac{\lambda + \nu -1}{i} \ffac{\mu_{+}}{\nu - i} \ffac{\kappa_{+} + \mu_{+} -1}{\nu - i} }{\rfac{\kappa_{-}}{i}} \, 
              \rho_x^{\mu_{+} - \nu +i} \rho_y^i \psi_{\kappa_{-} +i}
     \end{equation}
     for $\nu=0,1,2,\dots$.
\end{list}
\end{lemma}
\begin{proof}
We only show \eqref{e:iteration_X^{+}} by induction on $\nu$ here.
The other case \eqref{e:iteration_X^{-}} can be shown similarly.

It is trivial if $\nu=0$, and it is nothing but Propostion \ref{p:sl2-action_on_ourE} if $\nu=1$. 
Assume that it is true for $\nu \geqsl 1$, and apply $X^{+}$ to the both sides of \eqref{e:iteration_X^{+}}.
Then, one sees that the right-hand side equals 
\begin{equation}
\label{e:X^{+}_(nu+1)-times}  
\begin{aligned}
  & (-1)^{\nu} h_1 h_2 \sum_{i=0}^{\nu} \binom{\nu}{i} 
                \frac{ \ffac{\lambda+\nu-1}{i} \ffac{\mu_{+}}{\nu-i} \ffac{\kappa_{+}+\mu_{+}-1}{\nu-i} }{ \rfac{\kappa_{-}}{i} }
     \\
  & \times \left(%
              - (\mu_{+}-\nu+i) (\kappa_{+} + \mu_{+} - \nu + i - 1) \rho_x^{\mu_{+} - \nu + i - 1} \rho_y^i \psi_{\kappa_{-}+i}  
              + \frac{ \lambda+i-2\nu}{\kappa_{-}+i} \rho_x^{\mu_{+}-\nu+i} \rho_y^{i+1} \psi_{\kappa_{-}+i+1}
    \right).
\end{aligned}
\end{equation}
%
The coefficient of $(-1)^{\nu} h_1 h_2 \rho_x^{\mu_{+}-\nu+j-1} \rho_y^j \psi_{\kappa_{-}+j}$ in \eqref{e:X^{+}_(nu+1)-times}, $j=0,1,\dots,\nu+1$, equals
\begin{align*}
   & \binom{\nu}{j} \; \frac{ \ffac{\lambda+\nu-1}{j} \ffac{\mu_{+}}{\nu-j} \ffac{\kappa_{+}+\mu_{+}-1}{\nu-j} }{ \rfac{\kappa_{-}}{j} } 
         \cdot (-1) (\mu_{+}-n+j)(\kappa_{+}\mu_{+}-\nu+j-1)
          \\
   & 
     + \binom{\nu}{j-1} \frac{ \ffac{\lambda+\nu-1}{j-1} \ffac{\mu_{+}}{\nu-j+1} \ffac{\kappa_{+}+\mu_{+}-1}{\nu-j+1} }{ \rfac{\kappa_{-}}{j-1} }  
        \cdot \frac{-(\lambda+2\nu-j+1)}{\kappa_{-}+j-1}
          \\[2pt]
 = & - \left\{ \, \binom{\nu}{j} \, (\lambda+\nu-j) + \binom{\nu}{j-1} \, (\lambda+2\nu-j+1) \, \right\}
          \\
   & \hspace{5em} 
        \times  \frac{ \ffac{\mu_{+}}{\nu-j+1} \ffac{\kappa_{+}+\mu_{+}-1}{\nu-j+1} \ffac{\lambda+\nu-1}{j-1} }{ \rfac{\kappa_{-}}{j} }
          \\[2pt]
 = & \; - \binom{\nu+1}{j} \frac{ \ffac{\lambda+\nu}{j} \ffac{\mu_{+}}{\nu-j+1} \ffac{\kappa_{+}+\mu_{+}-1}{\nu-j+1} }{ \rfac{\kappa_{-}}{j} }.
\qedhere
\end{align*} 
\end{proof}

The following is our main result.

%
%
\begin{theorem}
\label{t:main_result}
Assume that $p \geqsl 1$, $q \geqsl 1$ and $p+q \in 2 \N$.
Let $m \in \N$ be a non-negative integer satisfying $m+3 \leqsl (p+q)/2$.
Then one has the following.
\begin{list}%
  {\upshape (\arabic{mylist})}{\setlength{\itemindent}{3pt} \setlength{\leftmargin}{3pt} \usecounter{mylist}}
    \item
\label{thm_item:K-type}

     The $K$-type formula of $\Module{\pm}{m}$ is given by 
     \begin{equation}
\label{e:K-type_Module_m}
     \left. \Module{\pm}{m} \right|_{K} 
       \simeq \bigoplus_{\substack{k, l \geqsl 0 \\  k - l + \frac{p-q}2  \in \Lambda_m}} \Harm^k(\R^p) \otimes \Harm^l(\R^q),
     \end{equation}  
     where $\Lambda_m=\{-m, -m+2, -m+4, \dots, m-2, m\}$, the set of $H$-weights of $F_m$;

    \item
\label{thm_item:irreducibility}
     Suppose further that $p, q \geqsl 2$.
     Then $\Module{\pm}{m}$ are irreducible $(\g,K)$-modules.
\end{list}
\end{theorem}
\begin{proof}
It suffices to show the theorem for $\Module{+}{m}$.
%
%
Let $f = h_1 h_2 \rho_y^{\mu_-} \psi_{\kappa_+} \ne 0$ be an element of $\Module{+}{m}$, where $h_1 \in \Harm^k(\R^p), h_2 \in \Harm^l(\R^q)$.
Then by Lemma \ref{l:iteration_sl2-action}, one obtains
\begin{align*}
  (X^-)^{m+1} f 
  &= h_1 h_2 \sum_{i=0}^{m+1} \binom{m+1}{i} \frac{ \ffac{0}{i} \ffac{\mu_{-}}{m+1 - i} \ffac{\kappa_{-} + \mu_{-} -1}{m+1 - i} }{ \rfac{\kappa_{+}}{i} } \,
           \rho_x^i \rho_y^{\mu_{-} - m-1 +i} \psi_{\kappa_{+} + i}
     \\
  &= \ffac{\mu_{-}}{m+1} \ffac{\kappa_{-} + \mu_{-} -1}{m+1} h_1 h_2 \rho_y^{\mu_{-} -m -1} \psi_{\kappa_{+}}.
\end{align*}
Thus, $(X^-)^{m+1} f = 0$ implies that $\ffac{\mu_{-}}{m+1} = 0$ or $\ffac{\kappa_{-}+\mu_{-}-1}{m+1} = 0$.
Namely,
\begin{align}
  \mu_{-} &= 0,1,\dots,m, \quad \text{or} \quad
\label{e:condi_on_mu_1}
     \\
  \mu_{-} &= -\kappa_{-}+1,\, -\kappa_{-}+2,\, \dots,\, -\kappa_{-}+m+1.
\label{e:condi_on_mu_2}
\end{align}
The assumption that $m+3 \leqsl (p+q)/2$, however, implies that \eqref{e:condi_on_mu_2} is impossible; 
if it holded true, then it would follow from \eqref{e:h_weight} that
\begin{align}
  m &=  -\kappa_{+} + \kappa_{-} + 2 \mu_{-}  
\label{e:m_kappa_plus_kappa_minus}
        \\
    &= -\kappa_{+} + \kappa_{-} + 2(-\kappa_{-} + i) 
\notag
        \\
    &= -(\kappa_{+}+\kappa_{-}) + 2 i
\notag
\end{align}
for $i=1,2,\dots,m+1$, hence that
\[
  \kappa_{+} + \kappa_{-} = -m + 2 i \leqsl m+2,
\] 
which contradicts $\kappa_{+} + \kappa_{-} \geqsl (p+q)/2 \geqsl m+3$.
Therefore, it follows from \eqref{e:condi_on_mu_1} that
\[
  k - l + \frac{p-q}{2} = \kappa_{+} - \kappa_{-} = -m + 2 \mu_{-} \in \Lambda_m,
\]
which proves (\ref{thm_item:K-type}).

%
%
Let us consider a closed subset $D_m \subset \R^2$ (with respect to the standard topology of $\R^2$) given by
\begin{equation}
\label{e:ambient_support_K-type_Module_m}
  D_m = \left\{ \left. (t_1,t_2) \in \R^2 \, \right| \,
                   t_1 \geqsl p/2,\,  t_2 \geqsl q/2, \, \abs{t_1 - t_2} \leqsl m 
        \right\},
\end{equation}
and the set of integral points of $D_m$ given by
\begin{equation}
\label{e:support_K-type_Module_m}
 \suppKtypeModule_m 
     = \left\{ (t_1,t_2) \in D_m 
                \left| 
                  \begin{array}{l}
                     t_1-p/2 \in \N, \, t_2-q/2 \in \N, \\
                     t_1 - t_2 \in \Lambda_m
	          \end{array}
                \right. 
        \right\}.
\end{equation}
Note that the sum in the right-hand side of \eqref{e:K-type_Module_m} can be written as the one with $(\kappa_{+},\kappa_{-})$ running over the set $\suppKtypeModule_m$.

Now, applying \eqref{e:p-action_on_HWV} to $f=h_1 h_2 \rho_y^{\mu_-} \psi_{\kappa_+} \in \Module{+}{m}$,
we denote the coefficient of 
\[
  (\pd_{x_i} h_1) (\pd_{y_j} h_2),\quad  (x_i h_1)(\pd_{y_j} h_2)^\dag,\quad  (\pd_{x_i} h_1) (y_j h_2)^\dag \quad \text{and} \quad (x_i h_1)^\dag (y_j h_2)^\dag
\]
in the right-hand side of \eqref{e:p-action_on_HWV} by $C_{--}$, $C_{+-}$, $C_{-+}$ and $C_{++}$ respectively, where $\mu_- = 0,1,\dots,m$.
Namely, 
\begin{align*}
  C_{--} &= \frac{\kappa_{-} + \mu_{-} - 1}{\kappa_{-} - 1} \rho_y^{\mu_{-}} \psi_{\kappa_{+} - 1},  
    &  
  C_{-+} &= \mu_{-} \, \rho_y^{\mu_{-} - 1} \psi_{\kappa_{+} - 1},
       \\
  C_{+-} &= \frac{\kappa_{+} - \mu_{-} - \kappa_{-}}{\kappa_{+} (\kappa_{-} - 1)} \rho_y^{\mu_{-}+1} \psi_{\kappa_{+} + 1},  
    & 
  C_{++} &= \frac{\kappa_{+} - \mu_{-} - 1}{\kappa_{+}} \rho_y^{\mu_{-}} \psi_{\kappa_{+} + 1}.
\end{align*}
(i) First, let us consider the case where $(\kappa_+,\kappa_-) \in \suppKtypeModule_m$ is an interior point of $D_m$.
Then, one obtains $\mu_- =1,2,\dots,m-1$ by \eqref{e:condi_on_mu_1} and \eqref{e:m_kappa_plus_kappa_minus}.
In particular, $C_{-+} \ne 0$. 
Now, $C_{--} = 0$ would imply $\mu_- = -\kappa_{-} + 1$, which contradicts $m+3 \leqsl (p+q)/2$ as we saw above. 
It also follows from \eqref{e:m_kappa_plus_kappa_minus} that $\kappa_{-} - \kappa_{+} + \mu_{-} = m - \mu_{-}$, and $C_{+-} \ne 0$. 
Finally, $C_{++} = 0$ would imply that $\kappa_{+} + \kappa_{-} = m+2$, which is absurd. 
Thus, all the coefficients in \eqref{e:p-action_on_HWV} never vanish.

\noindent
(ii) Next, let us consider the case where $(\kappa_+,\kappa_-) \in \suppKtypeModule_m$ is in the boundary of $D_m$.
Then there are three sub-cases:
\begin{enumerate}
\item[(ii-a)] $\mu_- = 0$,

\item[(ii-b)] $\mu_- = m$,

\item[(ii-c)] $0 < \mu_- < m$ and $k=0$ or $l=0$.
\end{enumerate}
In Case (ii-a), $C_{-+}=0$, and by the same reason as Case (i), $C_{--}$, $C_{+-}$ and $C_{++}$ are non-zero.
In Case (ii-b), $C_{+-}=0$ since $\kappa_{+}-\kappa_{-}=m$, and all the other coefficients are non-zero.
In Case (ii-c), all the coefficients are non-zero, but $\pd_{x_i} h_1 = 0$ or $\pd_{y_j} h_2 = 0$.

Therefore, one can move from any $K$-type in $\Module{+}{m}$ to any other $K$-type in $\Module{+}{m}$ by applying $\pi(X_{i,j}^-)$.
This completes the proof of (\ref{thm_item:irreducibility}), and of the theorem.
\end{proof}

\begin{example}
\label{ex:p=14,q=12,m=4}
Figure \ref{fig:suppKtypeModule} below illustrates $D_m$ in \eqref{e:ambient_support_K-type_Module_m} and $\Sigma_m$ in \eqref{e:support_K-type_Module_m} 
in the case where $p=14$, $q=12$ and $m=4$.
The colored area and the dots sitting in the area indicates $D_m$ and $\Sigma_m$ respectively.
Each $K$-type of $\Module{\pm}{m}$ corresponds to a dot by the correspondence
\[
  \Harm^k(\R^p) \otimes \Harm^l(\R^q) \longleftrightarrow (k,l) \quad (\text{or} \; (\kappa_{+},\kappa_{-})).
\]
Let us apply $\pi(X_{i,j}^{-})$ to an element $f$ of $\Module{\pm}{m}$. 
Then, if the $K$-type of $f$ corresponds to a white dot $\circ$ in Fig.~\ref{fig:suppKtypeModule}, 
one can move to any adjacent dots in the north-east, north-west, south-east, and south-west direction;
if it corresponds to a black dot $\bullet$ in Fig.~\ref{fig:suppKtypeModule}, 
one can move only to adjacent dots in the interior or in the boundary of $D_m$.

\begin{figure}[hbt]
\begin{tikzpicture}[scale=0.75, x=5mm,y=5mm, >=stealth] 

\path [fill=gray!35] (7,6) -- (7,11) -- (13,17) -- (17,13) -- (10,6) -- (7,6);

\draw[->] (0,0) -- (17,0) node[right] {$\kappa_{+}$};
\draw[->] (0,0) -- (0,17) node[above] {$\kappa_{-}$};
\draw[thick,->] (7,6) -- (17,6) node[right] {$k$};
\draw[thick,->] (7,6) -- (7,17) node[above] {$l$};

\draw[dashed] (7,6) -- (0,6) node[left] {\small $q/2$};
\draw[dashed] (7,6) -- (7,0) node[below] {\small $p/2$};

\draw (0,4) node[left] {$m$} -- (13,17) node[above right] {\small $\kappa_+ - \kappa_- = -m$};
\draw (0,2) -- (14,16) node[above right, rotate=-10] {\small $\ddots$}; 
\draw (0,0) node[below left] {$0$} -- (15,15) node[above right] {\small $\kappa_+ - \kappa_- = 0$};
\draw (2,0) -- (16,14) node[above right,rotate=-10] {\small $\ddots$}; 
\draw (4,0) node[below] {$m$} -- (17,13) node[above right] {\small $\kappa_+ - \kappa_- = m$};

\foreach \x in {7,...,12} \draw[fill] (\x,\x+4) circle (2pt); 

\draw[fill] (7,9) circle (2pt); 
\foreach \x in {8,...,13} \draw[fill=white] (\x,\x+2) circle (2pt);

\draw[fill] (7,7) circle (2pt); 
\foreach \x in {8,...,14} \draw[fill=white] (\x,\x) circle (2pt);

\draw[fill] (8,6) circle (2pt); 
\foreach \x in {9,...,15} \draw[fill=white] (\x,\x-2) circle (2pt);

\foreach \x in {10,...,16} \draw[fill] (\x,\x-4) circle (2pt); 


\end{tikzpicture}
\caption{Applying $\pi(X_{i,j}^{-})$, one can move from $\circ$ to dots in NE, NW, SE and SW directions, while from $\bullet$, only to dots in the interior or in the boundary.}
\label{fig:suppKtypeModule}
\end{figure}
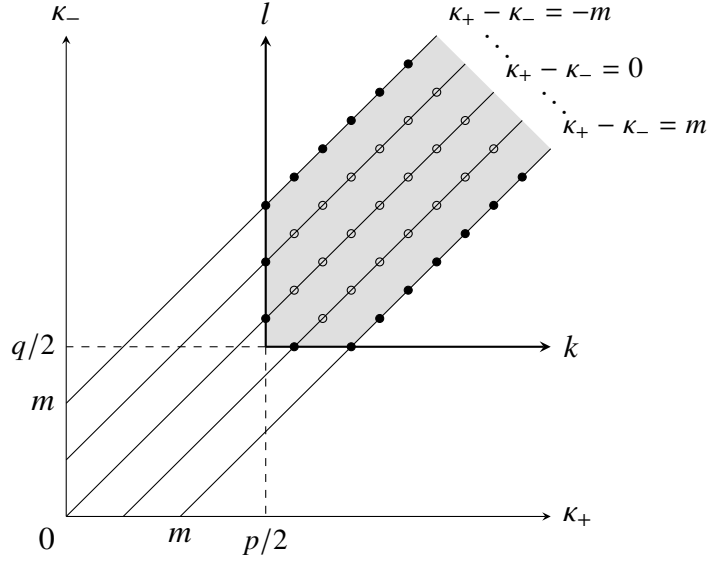

\end{example}

Now, let us briefly recall the definitions of the Gelfand-Kirillov dimension and the Bernstein degree of a finitely generated $U(\g)$-module $M$,
where $U(\g)$ denotes the universal enveloping algebra of $\g$.
Namely, we choose a finite-dimensional subspace $M_0$ so that $M=U(\g) M_0$, and for each non-negative integer $n$, we set $M_n := U_n(\g) M_0$
with $U_n(\g)$ denoting the subspace of $U(\g)$ spanned by products of at most $n$ elements of $\g$.
Then there exists a polynomial $\psi_{M}(t) \in \Q[t]$ of degree $d-1$ such that 
\begin{equation}
  \psi_{M}(n) = \dim (M_n/M_{n-1})  \quad \text{for all sufficiently large } n.
\notag
\end{equation}
Moreover, the leading term of $\psi_{M}$ is of the form
\begin{equation}
  \frac{m}{(d-1)!} t^{d-1}
\notag
\end{equation}
for a positive integer $m$.
We call $d$ the Gelfand-Kirillov dimension of $M$, and $m$ its Bernstein degree, which we denote by $\GKdim M$ and $\B-deg M$ respectively (see \cite{Vogan78} for more details).

%
%
\begin{corollary}
\label{c:GKdim_and_B-deg}
If $p,q \geqsl 2, p+q \in 2 \N$ 
and $m+3 \leqsl (p+q)/2$, 
then the Gelfand-Kirillov dimension and the Bernstein degree of \,$\Module{\pm}{m}$ are given by
\begin{align}
  \GKdim \Module{\pm}{m} &= p+q-3,
\label{e:GKdim}  
         \\
  \B-deg \Module{\pm}{m} &= \frac{ 4 (m+1) (p+q-4)!}{(p-2)! (q-2)!} 
\label{e:B-deg}
\end{align}
respectively.
\end{corollary}
\begin{proof}
Without loss of generality, one can assume that $p \geqsl q$.
We will consider $\Module{+}{m}$ here.
Then, let $\ell(j)$ be a line in $\R^2$ given by
\[
  \ell(j)=\left\{ (t_1, t_2) \in \R^2 \left| \, t_1 + t_2 = j \right. \right\} 
\]
with $j \in \N$, and set
\begin{equation}
\label{e:minimum_j}
  c:= \min \left\{ j \in \N \left| \, \ell(j) \cap \suppKtypeModule_m = \Lambda_m \right. \right\}. 
\end{equation}
As a generating ($K$-invariant) subspace of $\Module{+}{m}$,
we take
\begin{equation}
\label{e:M_0}
 \mathsf{M}_0 
     := \bigoplus_{ \substack{ (\kappa_{+}, \kappa_{-}) \in \suppKtypeModule_m \\ \kappa_{+}+\kappa_{-} \leqsl c} } 
          \Harm^k(\R^p) \otimes \Harm^l(\R^q) \rho_y^{\mu_-} \psi_{\kappa_+},
\end{equation}
where, in each summand, $\mu_{-}$ is determined by $\mu_{-}=(1/2)( \kappa_{+} - \kappa_{-} + m )$.
If one sets $\mathsf{M}_n := U_n(\g) \mathsf{M}_0$ ($\mathsf{M}_{-1}:=0$), then it follows from \eqref{e:K-type_Module_m} and \eqref{e:dim_Harm} that
\allowdisplaybreaks{
\begin{align}
\label{e:dim_M_n/M_{n-1}}
    & \dim (\mathsf{M}_n / \mathsf{M}_{n-1})  \notag \\
  = & \sum_{j=0}^m \dim \,( \Harm^{n+j}(\R^p) \otimes \Harm^{n+m-j+\frac{p-q}{2}}(\R^q) )
         \notag \\
  = & \, 4 \sum_{j=0}^m \frac{n+j+\tfrac{p}{2}-1}{(p-2)!} (n+j+1)(n+j+2) \cdots (n+j+p-3)
         \notag \\
    &  \hspace{1em}  \times  \frac{n+m-j+\tfrac{p-q}{2}+\tfrac{q}{2}-1}{(q-2)!} (n+m-j+\tfrac{p-q}{2}+1)(n+m-j+\tfrac{p-q}{2}+2) 
         \notag \\
    & \hspace{20em} \cdots (n+m-j+\tfrac{p-q}{2}+q-3)
         \notag \\
  = & \, \frac{4(m+1)}{(p-2)! (q-2)!} \, n^{p+q-4} + (\text{lower order terms in } n)
\end{align}
}%
for all $n \in \N$, which implies \eqref{e:GKdim}.
Furthermore, since the leading term of \eqref{e:dim_M_n/M_{n-1}} can be rewritten as
\[
  \frac{4(m+1)}{(p-2)! (q-2)!} \, n^{p+q-4} = \frac{4(m+1) (p+q-4)!}{(p-2)! (q-2)!} \frac{n^{p+q-4}}{(p+q-4)!},
\]
one obtains \eqref{e:B-deg}.
This completes the proof.
\end{proof}

\begin{remark}%
\label{r:minimum_j_2}
One can show that the non-negative integer $c$ in \eqref{e:minimum_j} is in fact equal to
\(
  \max \{ m+p, m+q \}.
\)
\end{remark}

It is well known that the Gelfand-Kirillov dimension of the minimal representation of $\rmO(p,q)$ is equal to $p+q-3$ (cf.~\cite{Kobayashi_Orsted_03-1-3, ZH97}).
The $K$-type formula \eqref{e:K-type_Module_m} for $m=0$ in Theorem \ref{t:main_result} tells us that 
$\Module{+}{0}=\Module{-}{0}$ corresponds to the $(\g,K)$-module of the minimal representation of $\rmO(p,q)$.
However, as we have seen in Corollary \ref{c:GKdim_and_B-deg}, the Gelfand-Kirillov dimension of $\Module{\pm}{m}$ is equal to $p+q-3$ not only for $m=0$ 
but for {\em any} $m \in \N$ satisfying $m+3 \leqsl (p+q)/2$. 
The Bernstein degree can distinguish the minimal representation from the others.

\bibliographystyle{amsalpha}
\bibliography{rep,geom}

\end{document}